\documentclass[12pt,twoside]{article}

\usepackage{a4}
\usepackage{exscale}
\usepackage{fancyhdr}
\usepackage[centertags]{amsmath}
\usepackage{amsfonts}
\usepackage{amssymb}
\usepackage{theorem}
\usepackage{ulem}
\usepackage{isolatin1}
\usepackage{epsfig}

\renewcommand{\sectionmark}[1]
                    {\markboth{Kapitel \thesection\ #1}{}}
\renewcommand{\sectionmark}[1]
                 {\markright{} }

\newtheorem{thm}{Theorem}
\newtheorem{definition}{Definition}       
\newtheorem{lem}[thm]{Lemma}
\newtheorem{corollary}[thm]{Corollary}

{\theorembodyfont{\rmfamily}}

\numberwithin{equation}{section}

\def\parmod{\parskip=2pt plus1pt minus1pt}

\newenvironment{einr}{\parmod 
                      \begin{list}{}
                        {\setlength{\rightmargin}{0cm}
                         \setlength{\leftmargin}{0,75cm}
                         \setlength{\labelwidth}{0cm}
                         \setlength{\parsep}{1pt}
                         \setlength{\itemsep}{1pt}
                         \setlength{\topsep}{1pt}
                         \setlength{\partopsep}{0pt}
                         \setlength{\labelsep}{0cm}
                         \setlength{\listparindent}{0pt}
                         \setlength{\itemindent}{0pt}}
                      \item[] \ignorespaces}
                     {\unskip \end{list}}

\newenvironment{proof}
        {\pagebreak[2] \vspace{-1pt}{\bf Proof.}  }
        {\hfill $\blacksquare$ \vspace{2pt}}

\newenvironment{proofof}[1]
        {\pagebreak[2] \vspace{-1pt}{\bf Proof of #1.}  }
        {\hfill $\blacksquare$ \vspace{2pt}}

\def\Im{\mathop{{\rm Im}}}

\def\fa{\mathcal{F}}
\def\ga{\mathcal{G}}

\def\ma{\mathcal{M}}

\def\nat{{\rm I\! N}}

\def\co{{\mathbb C}}

\def\dk{{\mathbb D}}

\def\gl{\left\{}
\def\gr{\right\}}
\def\kl{\left(}
\def\kr{\right)}
\def\kj{\overline}

\def\limn{\lim_{n\to\infty}}

\def\In{\subseteq}

\def\abb{\longrightarrow}

\renewcommand{\rho}{\varrho}
\renewcommand{\phi}{\varphi}
\renewcommand{\epsilon}{\varepsilon}

\def\beq{\begin{equation}}
\def\eeq{\end{equation}}
\def\beqar{\begin{eqnarray}}
\def\eeqar{\end{eqnarray}}
\def\beqaro{\begin{eqnarray*}}
\def\eeqaro{\end{eqnarray*}}
\def\bsat{\begin{thm}}
\def\esat{\end{thm}}
\def\blem{\begin{lem}}
\def\elem{\end{lem}}
\def\bkor{\begin{corollary}}
\def\ekor{\end{corollary}}
\def\bdefin{\begin{definition}}
\def\edefin{\end{definition}}
\def\bbew{\begin{proof}}
\def\ebew{\end{proof}}
\def\bbewo{\begin{proofof}}
\def\ebewo{\end{proofof}}

\renewcommand{\rho}{\varrho}
\renewcommand{\phi}{\varphi}

\begin{document}

\thispagestyle{plain}

\fancyhead[CE]{J. Grahl and S. Nevo}
\fancyhead[CO]{Extension of a Normality Result by Carath\'eodory}
\fancyhead[LE,RO]{\thepage}
\fancyhead[LO,RE]{}
\fancyfoot[CE,CO]{}

\begin{center}
{\LARGE \bf An Extension of a \\[12pt] Normality Result by Carath\'eodory } \\[18pt]
{\LARGE \it - Draft - } \\[18pt]
{\Large \it by Jürgen Grahl and Shahar Nevo}
\end{center}

\renewcommand{\thefootnote}{}

\footnote{Part of this work was supported by the German Israeli
Foundation for Scientific Research and Development (No.~G 809$-$234.6/2003).}



\section{An extension of Carath\'eodory's normality result}

The probably most essential and fundamental result in the theory of normal
families is Montel's Theorem which says that a family of functions meromorphic
in a domain in $\co$ which omit three distinct fixed values in $\kj\co$ is
normal; Schiff \cite{schiff} calls it the "Fundamental Normality Test"
(FNT). There are two natural directions to generalize this result:  
\begin{itemize}
\item[(1)]
Instead of fixed exceptional values, one may consider exceptional values
depending on the respective function in the family under consideration. Of
course, in this context one can hope for normality only under additional
assumptions on these exceptional values: It does not suffice that the
exceptional values are distinct; they have to be kind of "uniformly
distinct". The respective version of the  FNT is due to Carath\'eodory
\cite[p.~202]{cara}. Here $\chi$ denotes the chordal and $\sigma$
the spherical metric on $\kj\co$. 

\bsat\label{carath}
Let $\fa$ be a family of meromorphic functions on a domain $D$. Suppose there
exists an $\varepsilon>0$ such that each $f\in\fa$ omits three distinct values
$a_f,b_f,c_f\in\kj\co$ satisfying  
$$\sigma(a_f,b_f)\cdot \sigma(a_f,c_f)\cdot \sigma(b_f,c_f)\ge \varepsilon.$$
Then $\fa$ is normal in $D$. 
\esat

\item[(2)]  
Instead of exceptional values one may consider exceptional functions with
disjoint graphs, i.e. omitting each other. 

The case of meromorphic exceptional functions is almost trivial: If $a,b,c$
are meromorphic functions on a domain $D$ omitting each other and if each
$f\in\fa$ omits $a$, $b$ and $c$, then we consider the family $\ga$ of the 
cross ratio functions 
$$z\mapsto \frac{f(z)-a(z)}{f(z)-b(z)}\cdot \frac{c(z)-b(z)}{c(z)-a(z)} $$
all of which omit the values 0, 1 and $\infty$ in $D$. By the FNT we obtain
the normality of $\ga$, hence of $\fa$. 

But an analogous normality result even holds for exceptional functions
$a,b,c:D\abb \kj\co$ which are merely continuous (w.r.t. the spherical metric
on $\kj\co$) and which have disjoint graphs as Bargmann et al. \cite{bargbonk}
have shown with the help of Ahlfors' theory of covering surfaces. 
 \end{itemize}

In the present paper we combine these two directions of generalization by
considering meromorphic exceptional functions which depend on the
respective function in the family under consideration. Our main result
is the following.

\bsat\label{mainresult}
Let $\fa$ be a family of functions meromorphic in $\dk$ and
$\varepsilon>0$. Assume that for each $f\in\fa$ there exist functions
$a_f,b_f,c_f$ meromorphic in $\dk$ or $\equiv\infty$ such that $f$
omits the functions $a_f,b_f,c_f$ in $\dk$ and 
$$\sigma(a_f(z),b_f(z))\cdot \sigma(a_f(z),c_f(z))\cdot
\sigma(b_f(z),c_f(z))\ge \varepsilon$$
for all $z\in\dk$. Then $\fa$ is a normal family. 
\esat

This result no longer holds for exceptional functions which are merely
continuous on $\kj\co$ as the following counterexample shows: 

{\bf Counterexample:} Let $f_n(z):=e^{nz}$, $a_n\equiv 0$, $b_n\equiv\infty$,
$c_n(z):=-e^{i n \Im(z)}$ for all $n$. Then each $f_n$ omits the
continuous functions $a_n$, $b_n$ and $c_n$ and 
$$\sigma(a_n(z),b_n(z))=\frac{\pi}{2}, \qquad
\sigma(a_n(z),c_n(z))=\sigma(b_n(z),c_n(z))=\frac{\pi}{4}$$
for all $n\in\nat$, $z\in\dk$, but $(f_n)_n$ is not normal in $\dk$.

\section{An extension of Zalcman's Lemma}

In the proof of Theorem \ref{mainresult} we require an extension of
Zalcman's well-known rescaling lemma \cite{zalc}. For $p\in\nat$ we
define the projections 
$$\pi_j:(\ma(\dk))^p\abb\ma(\dk) \qquad (1\le j\le p)$$ 
by
$$\pi_j(f_1,\dots,f_p):=f_j \qquad\mbox{ for } \qquad (f_1,\dots,f_p)\in(\ma(\dk))^p.$$ 

\blem{\bf (Extension of Zalcman's Lemma)} \label{zalcmanlemma}
Let $p$ be a natural number and $\fa\In (\ma(\dk))^p$. Assume that
there exists a $j_0\in \gl1,\dots,p\gr$ such that the family
$\pi_{j_0}(\fa)$ is not normal at $z=0$. Then there exist sequences
$(f_n)_n=((f_{1,n},\dots,f_{p,n}))_n\In\fa$, $(z_n)_n\in\dk$ and
$(\varrho_n)_n\in]0;1[$ such that $\limn z_n=0$, $\limn \varrho_n=0$
and such that the sequences $(g_{j,n})_n$ defined by 
$$g_{j,n}(\zeta):=f_{j,n}(z_n+\rho_n\zeta)$$
for all $j=1,\dots,p$ converge locally uniformly in $\co$ to functions
$g_j\in\ma(\co)\cup\gl\infty\gr$ such that at least one of the
functions $g_1,\dots,g_p$ is not constant.  
\elem

\bbew{}
By Marty's theorem there exist sequences $(f_n)_n\In \fa$ and $(z_n^*)_n\In \dk$
such that $\limn z_n^*=0$ and $\limn f_{j_0,n}^\#(z_n^*)=\infty$. We
define 
$$r_n:=\kl f_{j_0,n}^\#(z_n^*) \kr^{- \frac{1}{2}} + 2 \cdot
|z_n^*|.$$
Then we have $\limn r_n=0$ and $\frac{|z_n^*|}{r_n} \leq \frac{1}{2}$
for all $n$. Furthermore, we define 
$$M_n:=\max_{|z|\le r_n} \kl1-\frac{|z|^2}{r_n^2}\kr\cdot\kl
f_{1,n}^\#(z)+\dots+f_{p,n}^\#(z)\kr.$$
Then we have 
$$M_n=\kl1-\frac{|z_n|^2}{r_n^2}\kr\cdot\kl f_{1,n}^\#(z_n)+\dots+f_{p,n}^\#(z_n)\kr$$
for certain $z_n\in U_{r_n}(0)$. From
$$r_n \cdot M_n
\ge \kl1-\frac{|z_n^*|^2}{r_n^2}\kr\cdot \kl f_{j_0,n}^\#(z_n^*)\kr^{\frac{1}{2}}
\ge \frac{3}{4}\cdot\kl f_{j_0,n}^\#(z_n^*)\kr^{\frac{1}{2}}$$
we see $\limn r_n M_n=\infty$. We define
$$\varrho_n:=\frac{r_n^2-|z_n|^2}{r_n^2 M_n} \qquad\mbox{ and }\qquad
R_n:= \frac{r_n-|z_n|}{\rho_n}.$$ 
Then 
$$ \frac{\rho_n}{r_n-|z_n|}
= \frac{r_n+|z_n|}{r_n^2 M_n} 
\leq \frac{2}{r_n M_n}
\longrightarrow 0 \quad (n\to \infty).$$
So we have $\limn \rho_n=0$ and $\limn R_n=\infty$. The functions  
$$g_{j,n}(\zeta):= f_{j,n}(z_n+\rho_n \zeta)$$
are meromorphic in $U_{R_n}(0)$ and satisfy
\beqaro
g_{j,n}^\#(\zeta) 
&=& \rho_n\cdot f^\#_{j,n}(z_n+\rho_n\zeta)
\le \frac{\varrho_n M_n}{1-\frac{1}{r_n^2}\cdot |z_n+\rho_n\zeta|^2}\\
&=& \frac{r_n^2-|z_n|^2}{r_n^2-|z_n+\rho_n\zeta|^2}
\leq \frac{r_n-|z_n|}{r_n-|z_n|-\rho_n |\zeta|}
\leq \frac{1}{1-\frac{R}{R_n}}<2
\eeqaro
for $|\zeta|\le R<\frac{1}{2}\cdot R_n$. 
So by Marty's theorem each sequence $(g_{j,n})_n$ is normal in
$U_R(0)$ for every $R>0$. Therefore, we may assume that $(g_{j,n})_n$
converges locally uniformly in $\co$ to some
$g_j\in\ma(\co)\cup\gl\infty\gr$ for every $j=1,\dots,p$. From 
$$g_{1,n}^\#(0)+\dots+g_{p,n}^\#(0)
=\rho_n\cdot \kl f_{1,n}^\#(z_n)+\dots+f_{p,n}^\#(z_n)\kr
=\rho_n\cdot M_n\cdot \frac{r_n^2}{r_n^2-|z_n|^2}=1$$
we finally obtain 
$g_1^\#(0)+\dots+g_p^\#(0)=1$, i.e. $g_j^\#(0)>0$ for some $j\in\gl
1,\dots,p\gr$ (not necessarily $j=p$) which means that not all of the
$g_j$ are constant. 
\ebew

Also Pang's extension of Zalcman's Lemma can be generalized in a similar way. 

One major disadvantage of this Lemma is the fact that it does not give
any control which of the limit functions $g_j$ are non-constant. (Of 
course, it is trivial that only those $g_j$ can be non-constant for
which $\pi_j(\fa)$ is not normal.) Proving a stronger version of this
Lemma where one can prescribe which $g_j$ is non-constant would be a
giant leap for the theory of normal families. In many possible
applications (for example, to a famous conjecture of Cartan and
Eremenko \cite{erem1,erem2}) 
it would even suffice if one could exclude the case $g_j\equiv \infty$.
On the other hand, in general one cannot expect that one can construct
several nonconstant limit functions by simultaneous rescaling. The 
deeper reason for this is the fact that if $(f_n)_n$ and $(g_n)_n$ are
sequences in $\ma(\dk)$ which are not normal at the origin, then one
cannot conclude that there exists a sequence $(z_n)_n$ in $\dk$ with
$\limn z_n=0$ such that both $(f_n^\#(z_n))_n$ and $(f_n^\#(z_n))_n$ are
unbounded as the following counterexample illustrates. 

{\bf Counterexample:} Consider the functions $f_n(z):=nz+\sqrt{n}$ and
$g_n(z):=-nz+\sqrt{n}$. Then 
$$f_n^\#(z)=\frac{n}{1+n|1+\sqrt{n} z|^2}, \qquad
g_n^\#(z)=\frac{n}{1+n|1-\sqrt{n} z|^2},$$
so 
$$f_n^\#\kl -\frac{1}{\sqrt{n}}\kr=g_n^\#\kl \frac{1}{\sqrt{n}}\kr=n$$ 
for all $n$. Hence by Marty's theorem both sequences $(f_n)_n$ and
$(g_n)_n$ are not normal at the origin. On the other hand, 
$f_n^\#(z_n)\to\infty$ implies $1+\sqrt{n}z_n\to 0$ while
$g_n^\#(z_n)\to\infty$ implies $1-\sqrt{n}z_n\to 0$ for
$n\to\infty$. Obviously, it is impossible to satisfy both conditions
simultaneously. If 
$f_n(z_n+\rho_n \zeta)=nz_n+\rho_n n \zeta + \sqrt{n}$ (with
$z_n,\rho_n$ as in Zalcman's Lemma) tends to some nonconstant limit
function, then $g_n(z_n+\rho_n \zeta)=-f_n(z_n+\rho_n \zeta)+\sqrt{n}$
tends to $\infty$ (and vice versa).

\section{Proof of Theorem \ref{mainresult}}

We start with the corresponding result for functions meromorphic in
$\co$. 

\blem \label{picardtype}
Let $a,b,c\in\ma(\co)\cup\gl\infty\gr$ and $\varepsilon>0$. Assume
that 
\beq\label{conditiondistance}
\sigma(a(z),b(z))\cdot \sigma(a(z),c(z))\cdot \sigma(b(z),c(z))\ge
\varepsilon
\eeq
for all $z\in\co$. Then $a$, $b$ and $c$ are constant. 
\elem

\bbew{}
In view of $\sigma(z,w)\le \pi$ for all $z,w\in\kj\co$ we have 
$$\sigma(a(z),b(z))\ge \frac{\varepsilon}{\pi^2}$$
for all $z\in\co$. If $a\equiv\infty$, then we deduce that $b$ is
bounded, hence constant by Liouville's theorem. If
$a\not\equiv\infty$, then we can conclude that 
$$|a(z)-b(z)|\ge \sigma(a(z),b(z))\ge \frac{\varepsilon}{\pi^2}$$
for all $z\in\co$ which implies that $a-b$ is constant. Now from
$$\kl 1+|a(z)|^2\kr\cdot \kl 1+|b(z)|^2\kr
=\frac{|a(z)-b(z)|^2}{(\chi(a(z),b(z)))^2}
\le\frac{|a(z)-b(z)|^2}{\kl\frac{2}{\pi}\sigma(a(z),b(z))\kr^2}
\le \frac{\pi^6}{4\varepsilon^2}\cdot |a(z)-b(z)|^2$$
for all $z\in\co$ we see that $a$ and $b$ itself are bounded, hence
constant. 

So we have shown that $a$ and $b$ are constant. In the same way we
obtain that $c$ is constant as well.
\ebew

\blem \label{exceptfunctnormal}
Let $\ga\In (\ma(\dk))^3$ and $\varepsilon>0$. Assume that
$$\sigma(a(z),b(z))\cdot \sigma(a(z),c(z))\cdot \sigma(b(z),c(z))\ge \varepsilon$$
for all $(a,b,c)\in\ga$ and all $z\in\dk$. Then the families
$\pi_1(\ga),\pi_2(\ga)$ and $\pi_3(\ga)$ are normal in $\dk$. 
\elem

\bbew
If $\pi_j(\ga)$ was not normal, then by Lemma \ref{zalcmanlemma} one
could construct functions $a,b,c\in\ma(\co)\cup\gl\infty\gr$ which
satisfy (\ref{conditiondistance}) for all $z\in\co$ such that one of
the functions $a,b,c$ is not constant. This contradicts Lemma
\ref{picardtype}. 
\ebew

{\bf Proof of Theorem \ref{mainresult}:}
We assume that $\fa$ is not normal in $\dk$. W.l.o.g. we may assume
that $\fa$ is not normal at $z=0$. Then by Lemma \ref{zalcmanlemma}
there exist sequences $(f_n)_n\in \fa$,
$(a_n)_n,(b_n)_n,(c_n)_n\In\ma(\dk)\cup\gl\infty\gr$, $(z_n)_n\in\dk$
and $(\varrho_n)_n\in]0;1[$ such that $\limn z_n=0$, $\limn \varrho_n=0$, 
$f_n$ omits the functions $a_n,b_n,c_n$, 
$$\sigma(a_n(z),b_n(z))\cdot \sigma(a_n(z),c_n(z))\cdot \sigma(b_n(z),c_n(z))\ge
\varepsilon$$
for all $z\in\dk$ and all $n$ and such that the sequences $(g_{n})_n$,
$(A_n)_n$, $(B_n)_n$ and $(C_n)_n$ defined by 
$$g_n(\zeta):=f_n(z_n+\rho_n\zeta),$$
$$A_n(\zeta):=a_n(z_n+\rho_n\zeta), \qquad
B_n(\zeta):=b_n(z_n+\rho_n\zeta),\qquad
C_n(\zeta):=c_n(z_n+\rho_n\zeta)$$
converge locally uniformly in $\co$ to functions
$g,A,B,C\in\ma(\co)\cup\gl\infty\gr$, resp., not all of which are
constant. Now Lemma \ref{exceptfunctnormal} ensures that $(a_n)_n$,
$(b_n)_n$ and $(c_n)_n$ are normal. This forces $A$, $B$ and $C$ to be 
constant. Therefore, $g$ is not constant. By Hurwitz's theorem, $g$
omits the three distinct constants $A$, $B$ and $C$. This contradicts
Picard's theorem.

\vspace{12pt}
\parbox[t]{90mm}{\it
 J\"urgen Grahl \\
 University of W\"urzburg \\
 Department of Mathematics \\
 W\"urzburg\\
 Germany\\
 e-mail: grahl@mathematik.uni-wuerzburg.de}
\parbox[t]{80mm}{\it
Shahar Nevo \\
Bar-Ilan University\\
Department of Mathematics\\
Ramat-Gan 52900\\
Israel\\
e-mail: nevosh@macs.biu.ac.il}

\end{document}